\documentclass[12pt]{article}
\usepackage{amsfonts,euscript,enumerate}

\newtheorem{thm}{Theorem}[section]

\newtheorem{exam}{Example}
\newtheorem{quest}{Question}
\newtheorem{rem}[thm]{Remark}
\newtheorem{defn}[thm]{Definition}

\newtheorem{prop}[thm]{Proposition}
\newtheorem{lem}[thm]{Lemma}
\newtheorem{cor}[thm]{Corollary}
\def\eqref#1{(\ref{#1})}
\def\qed{~\vrule height8pt width 5pt depth -1pt\medskip}

\begin{document}
\bibliographystyle{plain}

\begin{center}
{\bf \large Noetherian algebras over algebraically closed fields}\\
\vspace{7 mm} \centerline{Jason P. Bell}
\centerline{Department of Mathematics} \centerline{Simon Fraser University}
\centerline{8888 University Drive}
\centerline{
Burnaby, BC,  V5A 1S6}
\centerline{CANADA} 
\centerline{\tt jpb@math.sfu.ca}
\end{center}
\begin{abstract} Let $k$ be an uncountable algebraically closed field and let $A$ be a countably generated left Noetherian $k$-algebra.  Then we show that $A\otimes_k K$ is left Noetherian for any field extension $K$ of $k$.  We conclude that all subfields of the quotient division algebra of a countably generated left
Noetherian domain over $k$ are finitely generated extensions of $k$.  We give examples which show that $A\otimes_k K$ need not remain left Noetherian if the hypotheses are weakened. 
\end{abstract} 
\section{Introduction} Let $k$ be a field. 
A left Noetherian $k$-algebra $A$ is called \emph{strongly left Noetherian} if for any commutative Noetherian $k$-algebra $C$, the algebra $A\otimes_k C
$ is again left Noetherian.    A weaker property is the \emph{stably left Noetherian} property.  A left Noetherian $k$-algebra is called \emph{stably left Noetherian} if $A\otimes_k K$ is left Noetherian for any field extension $K$ of $k$.  There are many examples of such rings: the Weyl algebras, group algebras of polycyclic by finite groups, and enveloping algebras of finite dimensional lie algebras are just a few 
\cite[Corollary 9.1.8]{MR}.  Strongly right Noetherian
 and stably right Noetherian algebras are defined analogously and, just as with
Noetherian rings, a ring that is both stably left and stably right Noetherian is
 called stably Noetherian.  We have the inclusions 
\[ 
\{  
{\rm Noetherian} 
\}  \supset  \{ 
{\rm stably ~Noetherian } 
\}  \supset  \{ {\rm strongly ~Noetherian} \} . 
\] 
In \S 4, we give examples which show that each of these containments is proper.  DeJong \cite{ASZ} showed that graded Noetherian algebras over algebraically closed fields are stably Noetherian.  More specifically, he proved the following result.
   
\begin{thm}\label{thm: DeJong} (DeJong \cite[Theorem 5.1, p. 605-606]{ASZ}) Let $k$ be an algebraically closed field, let $G$ be a group and let $A$ be a locally finite, left Noetherian $G$-graded $k$-algebra.  Then $A$ is stably left Noetherian. 
\end{thm}  
DeJong's proof is geometric and relies heavily on the graded hypothesis.  We obtain a similar result in which we are able to avoid the graded hypothesis, but require instead the fact that the base field be uncountable in addition to being algebraically closed. 
Our main result is the following. 
\begin{thm} Let $A$ be a countably generated left Noetherian algebra over an uncountable algebraically closed field $k$.  Then $A$ is stably left Noetherian. \label{thm: main}
\end{thm} 
We obtain the following corollary from this result. 
\begin{cor} \label{cor: 1} Let $D$ be the quotient division algebra of a countably generated left Noetherian domain $A$ over an uncountable algebraically closed field $k$.  Then
all subfields of $D$ are finitely generated extensions of $k$. 
\end{cor}
In \S 2 we give some basic facts about the Zariski topology which we employ in the proof of the main theorem.   The facts we give are completely straightforward and one can safely skip this section, but we nevertheless include them for the sake of completeness.  
In \S 3, we prove Theorem \ref{thm: main}.  The main idea behind the theorem is
to assume that we have a countably generated but non-finitely generated left ideal $I$ in $A\otimes_k K$.  We pick a countable set of generators for $I$.  Then
the generators lie in a subalgebra of $A\otimes_k K$ of the form $A\otimes_k C$,
 where $C$ is a countably generated $k$-algebra.  We then use the fact that $C$
has many maximal ideals and every maximal ideal $M$ of $C$ satisfies $C/M\cong k
$.  For each maximal ideal $M$ of $C$, we look at the smallest number of generators needed to generate the image of $I$ in $A\otimes_k C/M\cong A$.  In this way
 we partition the maximal ideals of $C$ into countably many subsets.  We then argue that at least one subset must be dense and show that this is enough to deduce
that $I$ is finitely generated.   In \S 4, we show that the conclusion of the statement of Theorem \ref{thm: main} need not hold if any of the hypotheses are relaxed.

\section{Background on the Zariski topology}
In this section we give the basic facts about the Zariski topology that we will use in the proof of Theorem \ref{thm: main}.  Given a commutative ring $C$, we let ${\rm Spec}(C)$  denote the collection of prime ideals of $C$ and we let M-Spec$(C)$ denote the collection of maximal ideals of $C$.  We endow ${\rm Spec}(C)$ with a topology by declaring that a subset $Y$ of ${\rm Spec}(C)$ is closed if there is some ideal $I$ of $C$ such that $Y$ consists of all prime ideals of $C$ which contain $I$.  We then give M-Spec$(C)$ the subspace topology.  The main facts which we require are given in the following remarks.
\begin{rem} Let $C$ be a commutative algebra with Jacobson radical equal to $(0)$.  Then a subset $Y$ of ${\rm M}$-${\rm Spec}(C)$ is dense in ${\rm M}$-${\rm Spec}(C)$ if and only if
$\bigcap_{M\in Y} M = (0)$.\label{prop: Y}
\end{rem}
\vskip 2mm
Given a commutative ring $C$, we can think of an element $\alpha\in C$ as being a \emph{regular function} on M-Spec$(C)$, given by \begin{equation}
\alpha(M) = \alpha + M \ \in \ C/M. \end{equation}

\begin{rem} \label{prop: jac}
Let $C$ be a commutative algebra with Jacobson radical equal to $(0)$.  Suppose that
$\alpha\in C$ satisfies $\alpha(M)=0$ for all $M$ in some dense subset $Y$ of ${\rm M}$-${\rm Spec}(Y)$.  Then $\alpha=0$.
\end{rem}
\vskip 2mm
It is well-known that the Jacobson radical of a countably generated algebra over an uncountable field is nil (c.f. Amitsur \cite{A}); consequently, the Jacobson radical of a countably generated domain over an uncountable field is $(0)$.  We will also use this fact. 

\section{Algebras over uncountable fields}
\begin{defn} Let $X$ be a topological space.  We say that $X$ is \emph{countably irreducible} if $X$ cannot be expressed as a countable union of proper closed subsets.  Otherwise, we say that $X$ is \emph{countably reducible}.
\end{defn}
Amitsur \cite[Lemma 4, p. 41]{A} studied this property, showing that if $X$ is an affine irreducible variety over an algebraically closed, uncountable field $k$, then $X$ is countably irreducible.  Equivalently, if $k$ is an uncountable algebraically closed field and $A$ is a finitely generated domain over $k$, then it is impossible to construct a countable set of nonzero elements in $A$ such that every maximal ideal in $A$  contains at least one element from this set.   We extend Amitsur's result to countably generated algebras over uncountable fields.  We first use a well-known result.
\begin{lem} (Nullstellensatz) \label{lem: max} Let $k$ be an uncountable algebraically closed field and let $C$ be a countably generated commutative $k$-algebra.  Then every maximal ideal $M$ of $
C$ has the property $C/M\cong k$.  
\end{lem}
{\bf Proof.}  See McConnell and Robson \cite[Corollary 9.1.8]{MR}. \qed
\vskip 2mm
\begin{prop} \label{prop: 1} Let $k$ be an uncountable algebraically closed field and let $C$ be a countably generated commutative
$k$-algebra that is a domain.  Suppose $\theta_1,\theta_2,\ldots $ are nonzero elements of $C$.  Then there is a maximal ideal in $C$ not containing any of the $\theta_i$. 
\end{prop}
{\bf Proof.}   Let $S$ be the multiplicatively closed subset of $C$ generated by $\theta_1,\theta_2,\ldots $.  Then 
$B:=S^{-1}C$ is countably generated.  Pick a maximal ideal $M$ in $B$.  By Lemma \ref{lem: max}, $B/M\cong k$.   Note that we have an inclusion
$$C/(M\cap C) \hookrightarrow B/M \cong k.$$
Thus we have a map from $C/(M\cap C)$ into $k$.  This map must be surjective, since the image is a nonzero $k$-subspace of the image.  Hence $M'=M\cap C$ is a maximal ideal of $C$ which avoids $S$.  The result now follows. \qed
\vskip 2mm
An immediate consequence of Proposition \ref{prop: 1} is the following result.
\begin{cor} Let $k$ be an uncountable algebraically closed field and let $C$ be a countably generated commutative
$k$-algebra that is a domain.  If $C\not = k$, then ${\rm M}$-${\rm Spec}(C)$ is countably irreducible. \label{cor: count}
\end{cor}
{\bf Proof.} Suppose that M-Spec$(C)$ can be expressed as a countable union of proper closed sets.  Then there exist nonzero ideals $I_1,I_2, \ldots $ of $C$ such that every maximal ideal contains at least one of the $I_i$.  We pick nonzero $\theta_i\in I_i$.  Then every maximal ideal contains some $\theta_i$, contradicting Proposition \ref{prop: 1}.  \qed
\vskip 2mm
Let $k$ be an uncountable algebraically closed field and let $C$ be a countably generated commutative domain over $k$.  Given an element 
$\alpha\in C$ and a maximal ideal $M\subseteq C$, we note that by Lemma \ref{lem: max}, $C/M\cong k$.  Hence we can think of $\alpha$ as a regular function from ${\rm M}$-${\rm Spec}(C)$ to $k$ as follows; we define
\begin{equation} \alpha(M) := \alpha + M.\end{equation}
We now prove a key lemma.  This lemma shows that if there exist functions from the maximal spectrum of a commutative ring into a field $k$ which solve a system of linear equations defined over $C$ on a dense subset of the maximal spectrum, then there exists a solution to this system in the quotient field of $C$.
\begin{lem} \label{lem: sys} Let $k$ be an uncountable algebraically closed field, let $C$ be a countably generated $k$-algebra that is a domain and let $Y$ be a dense set of maximal ideals in  ${\rm M}$-${\rm Spec}(C)$.  Suppose there exist maps $u_1,\ldots ,u_d : Y\rightarrow k$ such that
$\sum_{i=1}^d \beta_{i,j}(M)u_i(M) \ = \ \alpha_i(M)$ for some elements
$\beta_{i,j},\alpha_i\in C$ with $1\le i\le d$ and $1\le j\le n$.
Then there exist elements $\gamma_1,\ldots ,\gamma_d$ in the quotient field of $C$ such that
$$\sum_{i=1}^d \beta_{i,j}\gamma_i \ = \ \alpha_i.$$
\end{lem}
{\bf Proof.}  Let $I$ denote the set of all polynomials $h(x_1,\ldots ,x_d)\in C[x_1,\ldots ,x_d]$ such that $h(u_1,\ldots ,u_d)(M)=0$ on $Y$.  We note that if $f(x_1,\ldots , x_d)$ and
$g(x_1,\ldots ,x_d)$ are elements of $I$ then $(f+g)(x_1,\ldots ,x_d)$ vanishes on $Y$.  Clearly $I$ is closed under multiplication by $C[x_1,\ldots ,x_d]$ and so it is an ideal.  We note also that $I\cap C=\{0\}$, for if  $\gamma\in I\cap C$, then $\gamma(M)=0$ for all $M$ in $Y$, which implies that $\gamma=0$, using Remark \ref{prop: jac} and the remarks which immediately follow its proof, since $Y$ is dense in ${\rm M}$-${\rm Spec}(C)$ and $C$ is a countably generated commutative domain.   We let $K$ denote the quotient field of $C$.  Then we can create the ring $B=K[x_1,\ldots ,K_d]/J$, where $J=KI$.  Let $\overline{x_i}$ denote the image of $x_i$ in $B$.  By construction, $B$ is reduced since if the square of a polynomial vanishes on $Y$ then it must vanish on $Y$ also.  Since
 \begin{equation}\Bigg(\sum_{i=1}^d \beta_{i,j}u_i -  \alpha_i\Bigg)(M) \ = \ 0 \qquad {\rm for~all~}M\in Y,
 \end{equation}
 and $Y$ is dense, Proposition \ref{prop: jac} gives
 \begin{equation}\sum_{i=1}^d \beta_{i,j}\overline{x_i} -  \alpha_i \ = \ 0.\end{equation}
 Since this relation holds in $B$, it holds in any homomorphic image of $B$.  Pick a maximal ideal $M$ of $B$.  Then by the Nullstellensatz \cite[Theorem 4.19, p. 132]{Ei}, $B/M$ is a finite extension of $K$.  Thus there exists a
 solution to the stystem of equations in unknowns $y_1,\ldots ,y_d$:
  \begin{equation} \sum_{i=1}^d \beta_{i,j}y_i -  \alpha_i \ = \ 0 \quad {\rm for} ~ j\le n.\end{equation}
   in a finite extension of $K$.  Since the $\beta_{i,j}$ and $\alpha_i$ all lie in $K$ and the system of equations is linear, we must have a solution in $K$.  The result now follows. \qed
   \vskip 2mm
   \noindent 
{\bf Proof of Theorem \ref{thm: main}.}
 Suppose $B=A\otimes_k K$ is not left Noetherian.  Then there exist elements 
 $b_1,b_2,\ldots \in B$ such that for each $i\ge 2$ we have
 \begin{equation}
 b_i \ \not \in \ \sum_{j<i} Bb_j.
 \end{equation}
Since $A$ is countably generated over $k$, we have that $${\rm dim}_k A \ \le \ \aleph_0.$$
Since finite dimensional algebras are stably Noetherian, $A$ is in fact countably infinite dimensional over $k$.  
Fix a $k$-basis \begin{equation}
 \mathcal{B} \ =\  \{r_1,r_2,\ldots \}\end{equation} for $A$.  Then
 for each $i$ we can write 
 \begin{equation}\label{eq: 1} b_i \ = \ \sum_{j=1}^{\infty} r_{j}\otimes \alpha_{i,j}\end{equation}
 for some elements
 $\alpha_{i,j}\in K$, where for $i$ fixed all but finitely many of the $\alpha_{i,j}$ are nonzero.
 We take $C$ to be the countably generated $k$-algebra generated by $\alpha_{i,j}$ with
 $i,j \ge 1$.  It is no loss of generality to assume that $K$ is the quotient field of $C$.  We note that the elements $b_1,\ldots $ can be regarded as elements of
 $A\otimes_k C$.   Given an element 
$r\in A\otimes_k C$ and a maximal ideal $M\subseteq C$, we note that by Lemma \ref{lem: max}, $C/M\cong k$.  Hence it makes sense to think of $r\in A\otimes_k C$ to be a map from ${\rm M}$-${\rm Spec}(C)$ to $A$ as follows; we define
 $r(M) \in A$ to be the image of $r$ under the natural map from $A\otimes_k C$ onto $A$ given by
$$A\otimes_k C \rightarrow A \otimes_k C/M \cong A.$$

Given a finite subset $S\subseteq \mathcal{B}$ and a natural number $i$, we let $X(i,S)$ denote the collection of maximal ideals $M$ of $C$ such that
$b_i(M) \ \in \ \sum_{j<i } V b_j(M),$
where $V$ is the $k$-vector space of $A$ spanned by the elements of $S$.  Given a maximal ideal 
$M$ of $C$, the left ideal $$Ab_1(M)+Ab_2(M)+\cdots $$ is finitely generated since $A$ is left Noetherian. Consequently, for every maximal ideal $M$, there exists some $i$ and some finite subset $S$ of $\mathcal{B}$ such that $M\in X(i,S)$.  We may assume that $C\not = k$ and thus M-Spec$(C)$ is countably irreducible by Corollary \ref{cor: count}.  Since there are only countably many finite subsets of $\mathcal{B}$, we conclude that there exist a finite subset $S$ of $\mathcal{B}$ and some natural number $i$ such that $X(i,S)$ is dense in M-Spec$(C)$.    By relabeling the elements of $\mathcal{B}$ if necessary, we may assume that $S=\{r_1,\ldots ,r_d\}$.  
Then we have maps $u_{\ell,j}:X(i,S)\rightarrow k$ for 
$\ell\le d$ and $j<i$ such that
\begin{equation}
b_i(M) \ = \ \sum_{\ell=1}^d \sum_{j<i} u_{\ell, j}(M)r_{\ell}b_j(M).
\end{equation}
(These maps $u_{\ell,j}$ are defined by choosing coefficients in $k$ which express $b_i(M)$ as a $k$-linear combination of $r_{\ell}b_j(M)$ for 
$1\le \ell\le d$ and $1\le j<i$.)
Note that 
equation (\ref{eq: 1}) gives
\begin{equation} \label{eq: 2}
\sum_{p=1}^{\infty} r_p \alpha_{i,p}(M) \ = \  \sum_{p=1}^{\infty} \sum_{j<i} \sum_{\ell =1}^d u_{\ell, j}(M)r_{\ell} r_{p}\alpha_{j,p}(M).
\end{equation}
For $\ell, p\ge 1$ we can express
\begin{equation} \label{eq: 3}
r_{\ell}r_p \ = \ \sum_{m=1}^{\infty} \beta_{\ell,p}^{(m)} r_m,\end{equation}
for some scalars $\beta_{\ell,p}^{(m)}\in k$.   Taking the coefficient of $r_m$ on both sides of equation (\ref{eq: 2}) using equation (\ref{eq: 3}) we see
\begin{equation} \label{eq: 4} \alpha_{i,m}(M) \ = \   \sum_{p=1}^{\infty} \sum_{j<i} \sum_{\ell =1}^d u_{\ell, j}(M)
\beta_{\ell, p}^{(m)} \alpha_{j,p}(M)
\end{equation}
for each $m\ge 1$.  We note that there exists some $N$ such that for every $j\le i$, $\alpha_{j,m}=0$ for all $m\ge N$.  Hence by Lemma \ref{lem: sys}, there exist elements $\gamma_{\ell,j}\in K$, $1\le \ell\le d$, $j<i$, such that
\begin{equation}  \label{eq: x}
\alpha_{i,m} \ = \   \sum_{p=1}^{\infty} \sum_{j<i} \sum_{\ell =1}^d \gamma_{\ell, j}
\beta_{\ell, p}^{(m)} \alpha_{j,p}, \end{equation}
for all $m\ge 1$.
Thus
\begin{eqnarray*}
b_i & = & \sum_{m=1}^{\infty} r_m\otimes \alpha_{i,m} \qquad {\rm~by~equation~(\ref{eq: 1})} \\
&=&  \sum_{m=1}^{\infty} r_m \otimes \Bigg( \sum_{p=1}^{\infty} \sum_{j<i} \sum_{\ell =1}^d 
\beta_{\ell, p}^{(m)} \otimes (\alpha_{j,p}\gamma_{\ell,j})\Bigg)\qquad {\rm~by~equation~(\ref{eq: x})}\\
&=& \sum_{p=1}^{\infty} \sum_{j<i} \sum_{\ell =1}^d 
\Bigg( \sum_{m=1}^{\infty} \beta_{\ell, p}^{(m)} r_m\Bigg)\otimes (\alpha_{j,p}\gamma_{\ell,j})\\
&=& \sum_{p=1}^{\infty} \sum_{j<i} \sum_{\ell =1}^d  r_{\ell}r_p\otimes (\alpha_{j,p}\gamma_{\ell,j})\qquad
{\rm by ~equation~(\ref{eq: 3})}\\
&=& \sum_{p=1}^{\infty} \sum_{j<i} \sum_{\ell =1}^d (r_{\ell}\otimes \gamma_{\ell,j})(r_p\otimes \alpha_{j,p})\\
&=&  \sum_{j<i} \sum_{\ell =1}^d (r_{\ell}\otimes \gamma_{\ell,j})\Bigg(  \sum_{p=1}^{\infty}r_p\otimes \alpha_{j,p}\Bigg)\\
&=& \sum_{\ell=1}^d \sum_{j<i} (r_{\ell} \otimes \gamma_{\ell,j}) b_j\qquad {\rm~by~equation~(\ref{eq: 1})}.
\end{eqnarray*}
But this contradicts the fact that $b_i\not \in Bb_1+\cdots + Bb_{i-1}$.   It follows that $A\otimes_k K$ is left Noetherian. \qed
\vskip 2mm
\noindent
{\bf Proof of Corollary \ref{cor: 1}.}  Let $K$ be a subfield of $D$.  By Theorem \ref{thm: main} $A\otimes_k K$ is left Noetherian.  Consequently, 
$D\otimes_k K$ is Noetherian since it is a localization of $A\otimes_k K$.  
We note that if $I\subset J$ are two ideals in 
$K\otimes K$ with $J$ properly containing $I$, then $(D\otimes_k K)I$ is properly contained by $(D\otimes_k K)J$ since $D\otimes_k K$
 is a free right $K\otimes_k K$-module.  Thus the lattice of left ideals of $K\otimes_k K$ embeds in the lattice of left ideals of $D\otimes_k K$ and so $K\otimes_k K$ is Noetherian.  A result of Vamos \cite{V} gives that $K$ is a finitely generated extension of $k$. \qed
\section{Examples} 
In this section we show that the conclusion of Theorem \ref{thm: main} need not
hold if any of the hypotheses are relaxed.  Furthermore, we show that the conclusion cannot be strengthened by replacing stably left Noetherian by strongly left
 Noetherian.  We accomplish this via a series of examples. 
\begin{exam} There exists an uncountable, non-algebraically closed field $k$ and a
 finitely generated left Noetherian $k$-algebra $A$ such that $A\otimes_k K$ is
not left Noetherian for some field extension $K$ of $k$. 
\end{exam} 
{\bf Proof.}  The following example is due to Resco and Small \cite{RS}.  Let $F
$ be an uncountable field of positive characteristic and let $K$ denote the
field extension $F(t_1,t_2,\ldots )$ of $F$ in countably many indeterminates $t_1,t_2, \ldots $.   Let $\delta: K \rightarrow K$ denote the $F$-derivation of $K
$ given by $\delta(t_i)=t_{i+1}$ for $i\ge 1$.  Let $A=K[x;\delta]$.  Then the centre of $A$ is the field 
$k=F(t_1^p,t_2^p,\ldots )$  
and $A$ is generated by  
$t_1$ and $x$ as a $k$-algebra.   
Furthermore, it is shown that $A\otimes_k K$ 
 is not left Noetherian. \qed \vskip 2mm 
\begin{exam} Let $k$ be an uncountable algebraically closed field.  Then there exists a stably Noetherian finitely generated $k$-algebra which is not strongly Noetherian. 
\end{exam} 
{\bf Proof.} Rogalski \cite{R} has shown that there exist elements $\alpha, \beta, \gamma \in k$ with $\alpha\beta =\gamma$ such that the subalgebra of $k\{x,y,z\}/(xy-\alpha yx,yz = \beta zy, xz = \gamma zx)$ generated by the images of $x-y$ and $y-z$ is Noetherian but not strongly Noetherian.  By Theorem \ref{thm: DeJong} this ring
is necessarily stably Noetherian. 
\qed \vskip 2mm
We note that any non-finitely generated field extension $K$ of a field $k$ is Noetherian but is not stably Noetherian since $K\otimes_k K$ is not left Noetherian \cite{V}.  Consequently, we see that the hypothesis that $A$ be countably generated is
necessary and that the hypothesis that $k$ be uncountable cannot be relaxed with
out also strengthening the hypothesis that $A$ be countably generated.  We thus pose the following question for finitely generated algebras. 
\begin{quest} Let $k$ be an algebraically closed field and let $A$ be a finitely
 generated left Noetherian $k$-algebra.  Is $A$ necessarily stably left Noetherian? 
\end{quest} 
We note that an affirmative answer to this question would immediately give Theorem \ref{thm: DeJong} of DeJong and hence could be viewed as a generalization of DeJong's result. 
\section{Acknowledgments} 
I thank John Farina, Dan Rogalski, and Lance Small for many interesting comments and suggestions.

\end{document}